\title{Periodic Orbits on Obtuse Edge Tessellating Polygons}
\author{Benjamin R. Baer, Faheem Gilani, Zhigang Han, and Ronald Umble}
\documentclass{article}
\usepackage{aliascnt}
\usepackage{fixltx2e,fix-cm}
\usepackage{amssymb}
\usepackage{amsmath}
\usepackage{amsfonts}
\usepackage{amsthm}
\usepackage{pgf,tikz}
\usepackage{mathrsfs}
\usepackage{epstopdf}
\usepackage{nameref,hyperref,cleveref,}
\usepackage{graphicx,caption,subcaption}
\usetikzlibrary{arrows}
\providecommand{\U}[1]{\protect\rule{.1in}{.1in}}
\theoremstyle{plain}

\def\NewTheorem#1
{	\newaliascnt{#1}{equation}
	\newtheorem{#1}[#1]{#1}
	\aliascntresetthe{#1}
	\expandafter\def\csname #1autorefname\endcsname{#1}
}
\NewTheorem{Theorem}
\NewTheorem{Proposition}
\NewTheorem{Corollary}
\NewTheorem{Lemma}
\numberwithin{equation}{section}

\newtheorem*{remark*}{Remark}
\begin{document}
\newpage
\maketitle

\section{Introduction}
Elementary mathematical billiards studies the motion of a massless particle moving with unit speed along a piece-wise linear path in the interior of a polygon $G$ subject to elastic reflections at the boundary $\partial G$, i.e., the angle of incidence equals the angle of reflection.  We think of $G$ as a frictionless billiard table, the edges of $G$ as its bumpers, the vertices of $G$ as its pockets, and the particle in motion as the cue ball. 

The particle's path is its \emph{orbit}. If the orbit reaches a vertex of $G$, it terminates and is \emph{singular}. A non-singular orbit that begins and ends at the same point is \emph{periodic} if the particle retraces its orbit when allowed to continue. A periodic orbit is \emph{primitive} if the particle traverses its orbit exactly once.  The \emph{period} of a periodic orbit is the number of times the particle strikes $\partial G$ as it traverses a primitive sub-orbit. 

In 2006, A. Baxter and R. Umble found and classified the periodic orbits on equilateral triangles \cite{Ba-Um}. Five year later, A. Baxter, E. McCarthy, and J. Eskreis-Winkler solved the analogous problem on rectangles, isosceles right triangles, and $30^\circ$-right triangles \cite{Ba-Es-Mc}. In this paper we solve the problem on $120^\circ$-isosceles triangles, $60^\circ$-rhombuses, and $60^\circ$-$90^\circ$-$120^\circ$-kites, and we make a conjecture in the case of regular hexagons.

An \emph{edge tessellation} of the plane is generated by reflecting a polygon $G$ and its reflected images in their edges. Each such $G$ lies in exactly one of the eight aforementioned families \cite{K-U}. The \emph{edges} of an edge tessellation are the edges of its polygons and the lines containing them are its \emph{inclines}. For example, an equilateral triangle in standard position generates an edge tessellation with inclines of $0^\circ$, $60^\circ$, and $120^\circ$. 

We identify a non-singular orbit in an edge tessellating polygon $G$ with a piece-wise linear curve $\gamma:I \rightarrow G$ defined on a finite interval $I$. An \emph{unfolding} of $\gamma$ in the edge tessellation $\mathcal{T}$ generated by $G$ is a line segment produced by successively reflecting $\gamma$ and its reflected images in the inclines of $\mathcal{T}$. Unfoldings relate the periodicity of $\gamma$ to the geometry of $\mathcal{T}$ and are sufficient for classifying periodic orbits in the non-obtuse cases. However, the analysis requires us to use more sophisticated techniques involving what we call the ``fence.''

\section{Periodic orbits on a $120^\circ$-isosceles triangle} \label{120isosceles}

\subsection{Preliminaries} 
Consider a $120^\circ$-isosceles $\triangle ABC$ positioned and labeled so that $\overline{AC}$ is horizontal with $A$ to the left of $C,$ the apex $B$ is positioned above $\overline{AC}$, and $\angle B$ is obtuse. The edge tessellation $\mathcal{T}$ generated by $\triangle ABC$ has inclines of $0^\circ$, $30^\circ,$ $60^\circ,$ $90^\circ,$ $120^\circ,$ and $150^\circ$. Given triangles $\triangle_1$ and $\triangle_2$ in $\mathcal T$ such that $\triangle_2=\tau (\triangle_1)$ for some translation $\tau$, two points $P_1$ on $\triangle_1$ and $P_2$ on $\triangle_2$ are \emph{(translationally) aligned (with respect to $\triangle_1$ and $\triangle_2$)} if $P_2=\tau (P_1)$. 

Let $\gamma:(0,T]\rightarrow \triangle ABC$ be an orbit in $\triangle ABC$. Let $\upsilon: (0,T]\rightarrow \mathcal{T}$ be an unfolding of $\gamma$ with \emph{initial point} $P = \lim_{t \searrow 0} \upsilon(t)$ and \emph{terminal point} $Q = \upsilon(T) = \sigma(\gamma(T))$, where $\sigma$ is the composition of reflections associated with $\upsilon$. The \emph{initial triangle} $\triangle ABC$ and the \emph{terminal triangle} $\triangle A'B'C' := \sigma \left( \triangle ABC \right)$ are \emph{consistently oriented} if $\sigma$ is orientation preserving. 

Unlike orientation, which is determined by comparing the labelings $(A,B,C)$ and $(A',B',C')=(\sigma(A),\sigma(B),\sigma(C))$ \cite{U-H}, alignment is determined by the relative positions of $P$ and $Q$ on $\triangle ABC$ and $\triangle A'B'C'$, and is independent of labeling. The periodicity of $\gamma$ is characterized by

\begin{Theorem}\label{align+orient}
 Let $\gamma$ be an orbit with unfolding $\upsilon$ whose initial point is on $\overline{AC}$. Then $\gamma$  is periodic if and only if the initial and terminal triangles are consistently oriented and the initial and terminal points are aligned.
\end{Theorem}

\begin{proof}
   When $\triangle A'B'C'$ and $\triangle ABC$ are consistently oriented, $\sigma$ is a rotation or a translation. Since $P$ and $Q$ are aligned, $\sigma$ is a translation and $\sigma (P)=Q$. Hence $P$ and $Q$ are in the same relative position in their respective triangles and the initial and terminal angles are equal. Therefore $\gamma$ is periodic.

 Conversely, given a periodic orbit $\gamma$, let $\upsilon$ be the unfolding with initial point $P$ on $\overline{AC}$, initial angle $\Theta : =m\angle QPC$, and  generated by first reflecting in a non-horizontal incline. If $\Theta\neq 90^\circ$, the reflection of $\overline{PQ}$ in the vertical line through $P$ is another unfolding, and we may restrict our considerations to initial angles in the range $0<\Theta\leq 90^\circ$. 

The images of $\overleftrightarrow{AC}$ in $\mathcal{T}$ are the horizontal, $60^{\circ}$, and $120^{\circ}$ inclines. We claim that $Q$ is on a horizontal incline. Suppose $Q$ is on a $60^\circ$ incline. Then the terminal angle at $Q$, which equals the initial angle $\Theta$ at $P$, is $60^\circ-\Theta$, $\Theta+120^\circ$, $\Theta-60^\circ$, or $240^\circ-\Theta$. But $0<\Theta\leq 90^\circ$ implies $\Theta=30^\circ$ (see Figure~\ref{initial-terminal-angle}), and an initial angle of $\Theta=30^\circ$ produces a period $8$ orbit that terminates on a horizontal incline, which is a contradiction (see Figure~\ref{30-unfolding}). Suppose $Q$ is on a $120^{\circ}$ incline.  Then the terminal angle at $Q$ is $120^\circ-\Theta$, $\Theta+60^\circ$, $\Theta-120^\circ$, or $300^\circ-\Theta$. But $0<\Theta\leq 90^\circ$ implies $\Theta=60^\circ$, and an initial angle of $\Theta=60^\circ$ produces an orbit of period $4$ or period $10$, both of which terminate on a horizontal incline, which is a contradiction (see Figure~\ref{alignex2}). Therefore $Q$ is on a horizontal incline as claimed.

\begin{figure}
\centering
\begin{minipage}[b]{.4\textwidth}
  \centering
  \definecolor{dots}{rgb}{0.8274509803921568,0.1843137254901961,0.1843137254901961}
\definecolor{niceblue}{rgb}{0,0.2,0.6}
\scalebox{0.9}{
\begin{tikzpicture}[line cap=round,line join=round,>=triangle 45,x=1cm,y=1cm]
\draw [shift={(4,0)},fill=black,fill opacity=0.05] (0,0) -- (0:0.38699481958058) arc (0:60:0.38699481958058) -- cycle;
\draw [shift={(1,0)},fill=black,fill opacity=0.05] (0,0) -- (0:0.38699481958058) arc (0:31.57482871319478:0.38699481958058) -- cycle;
\draw [shift={(5.65,2.857883832488647)},fill opacity=0.05] (0,0) -- (211.5748287131948:0.38699481958058) arc (211.5748287131948:240:0.38699481958058) -- cycle;
\draw [color=niceblue] (0,0)-- (5,0);
\draw [color=niceblue] (3.5,-0.8660254037844386)-- (6,3.4641016151377544);
\draw [,color=niceblue] (1,0)-- (5.65,2.857883832488647);
\draw (2,2) node[anchor= west] {Unfolding};
\draw (4.3,0.7) node[anchor=north west] {${\scriptstyle 60^{\circ}}$};
\draw (1.6370371521141402,0.49) node[anchor=north west] {${\scriptstyle \theta}$};
\draw [shift={(1,0)}] (0:0.38699481958058) arc (0:31.57482871319478:0.38699481958058);
\draw [shift={(1,0)}] (0:0.309595855664464) arc (0:31.57482871319478:0.309595855664464);
\draw (3.988004342791147,1.965074640549045) node[anchor=north west] {${\scriptstyle 60^{\circ}-\theta}$};
%\draw (4.088004342791147,2.165074640549045) node[anchor=north west] {${\scriptstyle 60^{\circ}-\theta}$};
\draw (0.7340492397594534,0.6041955349073744) node[anchor=north west] {$\mathit{P}$};
\draw (5.777881721626347,3.158361344139199) node[anchor=north west] {$\mathit{Q}$};
\draw (4.017002615984674,-0.3919909960021322) node[anchor=west] {Incline};
\draw [shift={(5.65,2.857883832488647)}] (211.5748287131948:0.38699481958058) arc (211.5748287131948:240:0.38699481958058);
\draw [shift={(5.65,2.857883832488647)}] (211.5748287131948:0.309595855664464) arc (211.5748287131948:240:0.309595855664464);
\draw [shift={(5.65,2.857883832488647)}] (211.5748287131948:0.232196891748348) arc (211.5748287131948:240:0.232196891748348);
\draw (5.494085520600588,0.3590988158396741) node[anchor=north west] {$\overleftrightarrow{\mathit{AC}}$};
\begin{scriptsize}
\draw [fill=dots] (1,0) circle (1.5pt);
\draw [fill=dots] (5.65,2.857883832488647) circle (1.5pt);
\end{scriptsize}
\end{tikzpicture}
}
    \caption{\small Equal initial and terminal angles of an unfolding.}
    \label{initial-terminal-angle}
\end{minipage}%
\hskip2cm\begin{minipage}[b]{.4\textwidth}
  \centering
 \definecolor{dots}{rgb}{0.8274509803921568,0.1843137254901961,0.1843137254901961}
\definecolor{niceblue}{rgb}{0,0.2,0.6}
%0.57735026919
\scalebox{1}{
\begin{tikzpicture}[line cap=round,line join=round,>=triangle 45,x=1cm,y=1cm]

\draw [color=niceblue] (0,0) -- (2,0);
\draw  [color=niceblue]  (2,0)-- (1,0.5773502691896257);
\draw  [color=niceblue]  (0,0)-- (1,0.5773502691896257);
\draw  [color=niceblue] (1,0.5773502691896257)-- (1,1.7320508075688772);
\draw  [color=niceblue]  (1,1.732050807568877)-- (2,0);
\draw  [color=niceblue] (2,1.1547005383792506)-- (1,1.7320508075688763);

\draw  [color=niceblue] (2,1.154700538379251)-- (2,0);
\draw  [color=niceblue] (3,0.57735026919)-- (2,0);
\draw  [color=niceblue]  (3,0.57735026919)-- (3,3*0.57735026919);
\draw  [color=niceblue]  (2,0)-- (3,3*0.57735026919);
\draw  [color=niceblue] (2,2*0.57735026919)-- (3,3*0.57735026919);
\draw  [color=niceblue]  (3,3*0.57735026919)-- (5,3*0.57735026919);
\draw  [color=niceblue]  (3,0.57735026919)-- (4,0);
\draw  [color=niceblue] (3,3*0.57735026919)-- (4,0);
\draw  [color=niceblue]  (3,3*0.57735026919)-- (4,2*0.57735026919);
\draw  [color=niceblue] (4,0)-- (4,2*0.57735026919);
\draw  [color=niceblue]  (5,3*0.57735026919)-- (4,2*0.57735026919);
\draw  [color=niceblue] (3,3*0.57735026919)-- (4,4*0.57735026919);
\draw  [color=niceblue]  (5,3*0.57735026919)-- (4,4*0.57735026919);

\draw  [color=niceblue]  (.7,0)-- (1+.7/2,.375278);
\draw [color=niceblue]  (1+.7/2,0)-- (1+.7/2,.375278);
\draw  [color=niceblue]  (.7,0)-- (.7/2,.35*.57735026919);
\draw [color=niceblue] (.35,0)-- (.7/2,.35*.57735026919);
\draw [dotted,color=niceblue] (1+.7/2,.375278)-- (3.7,3*.57735026919);
\begin{scriptsize}
\draw [fill=dots] (0.7,0) circle (1pt) node[anchor=north]{$P$};
\draw (0,0) node[anchor=north]{$A$};
\draw (2,0) node[anchor=north]{$C$};
\draw (1,0.5773502691896257) node[anchor=south east]{$B$};

\draw [fill=dots] (3.7,3*.57735026919) circle (1pt) node[anchor=north west]{$Q$};
\draw (3,3*.57735026919) node[anchor=south east]{$A'$};
\draw (5,3*.57735026919) node[anchor=south west]{$C'$};
\draw (4,4*0.5773502691896257) node[anchor=south]{$B'$};

\end{scriptsize}
\end{tikzpicture}
}
	\caption{\small An unfolding with initial angle $30^\circ$.}
	\label{30-unfolding}
\end{minipage}
\end{figure}

\begin{figure}[ht!]
	\centering
	
	\begin{subfigure}[b]{0.4\textwidth}
    \definecolor{dots}{rgb}{0.8274509803921568,0.1843137254901961,0.1843137254901961}
\definecolor{niceblue}{rgb}{0,0.2,0.6}
%0.57735026919
\scalebox{1}{
\begin{tikzpicture}[line cap=round,line join=round,>=triangle 45,x=1cm,y=1cm]

\draw [color=niceblue] (0,0) -- (2,0);
\draw  [color=niceblue]  (2,0)-- (1,0.5773502691896257);
\draw  [color=niceblue]  (0,0)-- (1,0.5773502691896257);
\draw  [color=niceblue] (1,0.5773502691896257)-- (1,1.7320508075688772);
\draw  [color=niceblue]  (1,1.732050807568877)-- (2,0);
\draw  [color=niceblue] (2,1.1547005383792506)-- (1,1.7320508075688763);

\draw  [color=niceblue] (2,1.154700538379251)-- (2,0);

\draw  [color=niceblue] (2,1.154700538379251)-- (3,3*0.57735026919);
\draw  [color=niceblue] (1,3*0.57735026919)-- (3,3*0.57735026919);
\draw  [color=niceblue] (1,3*0.57735026919)-- (2,4*0.57735026919);
\draw  [color=niceblue] (3,3*0.57735026919)-- (2,4*0.57735026919);
\draw  [color=niceblue]  (0,0)-- (2,6*0.5773502691896257);
\draw  [color=niceblue]  (2,4*0.5773502691896257)-- (2,6*0.5773502691896257)--(3,3*0.5773502691896257)--(3,5*0.5773502691896257)--(2,6*0.5773502691896257)--(4,6*0.5773502691896257)--(3,5*0.5773502691896257);

\draw  [color=niceblue]  (2,6*0.5773502691896257)-- (3,7*0.5773502691896257)--(4,6*0.5773502691896257);

\draw  [ color=niceblue] (.5,0)-- (.375, .216506);
\draw  [ color=niceblue] (.5,0)-- (3/4, .43301270189);
\draw  [ color=niceblue] (.75, .43301270189)-- (5/4, .43301270189);
\draw  [ color=niceblue] (1.5,0)-- (5/4, .43301270189);
\draw  [ color=niceblue]  (1.5,0)-- (2-.375,.216506);

\draw  [dotted, color=niceblue] (.75, .43301270189)-- (2.5,6*.57735026919);

\begin{scriptsize}
\draw [fill=dots] (0.5,0) circle (1pt) node[anchor=north]{$P$};
\draw (0,0) node[anchor=north]{$A$};
\draw (2,0) node[anchor=north]{$C$};
\draw (1,0.5773502691896257) node[anchor=south east]{$B$};

\draw [fill=dots] (2.5,6*.57735026919) circle (1pt) node[anchor=north west]{$Q$};
\draw  (2,6*0.5773502691896257) node[anchor=east]{$A'$};
\draw  (4,6*0.5773502691896257) node[anchor=west]{$C'$};
\draw   (3,7*0.5773502691896257) node[anchor=south]{$B'$};

\end{scriptsize}
\end{tikzpicture}
}
		\caption{An unfolding with initial angle $60^\circ$ and period 10.}
		\label{60-unfolding-10}
        \end{subfigure}
        \hskip1cm\begin{subfigure}[b]{0.4\textwidth}
		\definecolor{dots}{rgb}{0.8274509803921568,0.1843137254901961,0.1843137254901961}
\definecolor{niceblue}{rgb}{0,0.2,0.6}

%0.57735026919
\scalebox{1}{
\begin{tikzpicture}[line cap=round,line join=round,>=triangle 45,x=1cm,y=1cm]

\draw [color=niceblue] (0,0) -- (2,0);
\draw  [color=niceblue]  (2,0)-- (1,0.5773502691896257);
\draw  [color=niceblue]  (0,0)-- (1,0.5773502691896257);
\draw  [color=niceblue] (1,0.5773502691896257)-- (1,1.7320508075688772);
\draw  [color=niceblue]  (1,1.732050807568877)-- (2,0);
\draw  [color=niceblue] (2,1.1547005383792506)-- (1,1.7320508075688763);

\draw  [color=niceblue] (2,1.154700538379251)-- (2,0);

\draw  [color=niceblue] (2,1.154700538379251)-- (3,3*0.57735026919);
\draw  [color=niceblue] (1,3*0.57735026919)-- (3,3*0.57735026919);
\draw  [color=niceblue] (1,3*0.57735026919)-- (2,4*0.57735026919);

\draw  [color=niceblue] (3,3*0.57735026919)-- (2,4*0.57735026919);

\draw  [color=niceblue] (.8,0)-- (1.1, .519615242270663);
\draw  [dotted, color=niceblue] (1.8,3*0.57735026919)-- (1.1, .519615242270663);
\draw  [ color=niceblue] (.8,0)-- (.6, .346410);
\begin{scriptsize}
\draw [fill=dots] (0.8,0) circle (1pt) node[anchor=north]{$P$};
\draw (0,0) node[anchor=north]{$A$};
\draw (2,0) node[anchor=north]{$C$};
\draw (1,0.5773502691896257) node[anchor=south east]{$B$};

\draw [fill=dots] (1.8,3*.57735026919) circle (1pt) node[anchor=north west]{$Q$};
\draw  (1,3*0.57735026919) node[anchor=east]{$A'$};
\draw  (3,3*0.57735026919) node[anchor=west]{$C'$};
\draw  (2,4*0.57735026919) node[anchor=south]{$B'$};

\end{scriptsize}
\end{tikzpicture}
}
		\caption{An unfolding with initial angle $60^\circ$ and period 4.}
		\label{60-unfolding-4}
	\end{subfigure}      
	
	\caption{}
    \label{alignex2}
\end{figure}

Since $Q$ is on a horizontal incline, the base $\overline{A'C'}$ is horizontal. Furthermore, since $\overline{PQ}$ does not cross the interior of $\triangle A'B'C'$, its apex $B'$ lies above the base. Thus $\sigma$ is a translation or a glide reflection.  If $\sigma$ is a glide reflection, it reverses orientation. Then $180^\circ -\Theta =\Theta$ implies $\Theta=90^{\circ}$. But a periodic orbit with initial angle $90^\circ$ coincides with the period $8$ orbit with initial angle $30^{\circ}$. Since $\sigma$ is a composition of eight reflections, it preserves orientation, which is a contradiction. Therefore $\sigma$ is a translation, $P$ and $Q$ are aligned, and $\triangle A'B'C'$ and $\triangle ABC$ are consistently oriented.
\end{proof}

\begin{Corollary}\label{evenperiod}
The period of a periodic orbit is even, and an unfolding of a periodic orbit with initial point on a horizontal terminates on a horizontal.
\end{Corollary}

Corollary \ref{evenperiod} does not hold for all edge tessellating polygons. For example, Fagnano's periodic orbit on an equilateral triangle has period $3$ and terminates on a non-horizontal incline \cite{Ba-Um}.

Note that the periodic orbit displayed in Figure~\ref{30-unfolding}\text{ }  has initial angle $30^{\circ}$ and period $8$, while the periodic orbits displayed in Figure~\ref{alignex2}\text{ }  have initial angle $60^{\circ}$ and respective periods $4$ and $10$. A periodic orbit $\gamma$ is \emph{monoperiodic} if every periodic orbit with the same initial angle as $\gamma$ has the same period;  otherwise, $\gamma$ is \emph{biperiodic}. 

Our next proposition allows us to restrict our considerations to periodic orbits with initial angles
between 60 and 90 degrees.

\begin{Proposition}\label{align} 
Every periodic orbit can be represented by an unfolding with an initial angle $\Theta$  in the range $60^\circ \leq \Theta \leq 90^\circ$.
		
\end{Proposition}

\begin{proof}
Let $\gamma$ be a periodic orbit with initial angle $\Theta$ in the range $0^\circ<\Theta\leq 90^\circ$. Suppose $0^\circ<\Theta\leq 30^\circ$. Since $Q$ is on a horizontal incline, $\overline{PQ}$ cuts a $120^{\circ}$ incline at a point $P'$ with angle of incidence $\Phi=60^\circ+\Theta.$   Thus, there is an unfolding $\overline{P'Q'}$ of $\gamma$ with initial angle $\Phi$ in the range $60^\circ<\Phi\leq 90^\circ$. On the other hand, suppose $30^\circ<\Theta< 60^\circ$, and let $\overline{PQ'}$ be the reflection of $\overline{PQ}$ in $\overleftrightarrow{AC}.$ Then $\overline{PQ'}$ cuts a $60^{\circ}$ incline at a point $P'$ with angle of incidence $\Psi=120^\circ-\Theta$. Thus, there is an unfolding $\overline{P'R}$ with initial angle $\Psi$ in the range $60^\circ<\Psi<90^\circ$. 
\end{proof}

Since periodic orbits with initial angles $60^\circ$ and $90^\circ$ are understood, our problem reduces to classifying periodic orbits with initial angles in $(60^\circ, 90^\circ)$. 

\subsection{Contact points of an unfolding and the fence}

Let  $AC$ denote the length of $\overline{AC}$ and  let $u=\frac{1}{2}AC$. Impose a rectangular coordinate system on $\mathcal T$ with horizontal axis $\overleftrightarrow{AC}$, origin $O$ at the midpoint of $\overline{AC}$, horizontal unit of length $u$, and vertical unit of length $\sqrt{3}u$. Then points on vertical inclines have integer horizontal coordinates, points on horizontal inclines have integer vertical coordinates, and adjacent vertical and adjacent horizontal inclines lie one unit apart. If $\upsilon$ is an unfolding of a periodic orbit with initial point $P$, terminal point $Q$, and initial angle $\Theta$, Theorem~\ref{align+orient} implies that the vector $\mathbf{PQ}$ is parallel to a vector $\left(  x,y\right) $ for some $x,y\in \mathbb{N}$; hence $\Theta=\arctan\left(\frac{y}{x}\sqrt{3}\right)$. Furthermore, $\Theta\in\left(60^\circ, 90^\circ\right)$ implies $x < y$, and we may assume $\gcd\left(  x,y\right)  =1$. 

Parametrize $\upsilon$ via $\upsilon\left( t\right)  :=\left(  t+a,\frac{y}{x}t\right),$ $0< t \leq T$, where $a\in(-1,1)$; then $P = \lim_{t \searrow 0} \upsilon(t)=(a,0)$ and $Q = \upsilon(T)$. 
Each point at which $\upsilon$ cuts a vertical incline lies on a fundamental vertical segment of length 2 connecting the midpoints of two horizontal edges. The function $f:\mathbb{Z} \times\mathbb{R}\rightarrow \mathbb{R}/2\mathbb{Z}$
defined by $f\left(\alpha, \beta \right)  :=\alpha+ \beta +2 \mathbb Z$ identifies each such segment with the quotient group $\mathbb{R}/2\mathbb{Z}$, called the \emph{fence}. 
Geometrically, the projection in the $60^\circ$ direction onto the vertical coordinate axis sends $(\alpha,\beta) $ into the coset $f (\alpha,\beta) $. It will often be convenient to think of the fence as the interval $\mathcal{F}:=\left[0,2\right)$ of coset representatives and to write $f\left(\alpha, \beta \right) =\left(\alpha+ \beta \right) \operatorname{mod}2$. The fence $\mathcal{F}$ consists of the \emph{barrier} $\mathcal{B}:=\left(  \frac{1}{3},\frac{5}{3}\right]  $ and the \emph{gate} $\mathcal{B}^{c}:=\mathcal{F\smallsetminus B}$.

Between consecutive horizontal inclines, an unfolding $\upsilon$ always cuts four non-vertical edges. But whether or not $\upsilon$ also cuts a vertical edge is determined by its set of \emph{contact points} 
\begin{equation*}
\mathcal{C}_{T}:=\left\{  f\left(\upsilon(t)\right)  :  0< t \leq T \,\, \mbox{and} \,\, 
t+a \in \mathbb Z\right\}\subset \mathcal F.
\end{equation*}
Indeed, $\upsilon$ cuts a vertical edge at $\upsilon(t_i)$ if and only if $f(\upsilon(t_i))\in \mathcal B$. 
The \emph{multiplicity} of a contact point $c\in \mathcal{C}_{T}$, denoted by $m_{T}\left(  c\right) $, is the number  of times $\upsilon$ cuts a vertical incline at a point corresponding to $c$, i.e.,   
\begin{equation*}
m_{T}\left(  c\right)  :=\#\left\{  t :0< t \leq T \,\, \mbox{and} \,\, f\left(\upsilon\left(t \right)\right)  =c\right\}  .
\end{equation*}
These ideas are illustrated in Figure~\ref{figure4}. 

\begin{figure}[ht!]
	\centering
	
	\begin{subfigure}[b]{0.3\textwidth}
    \centering
    \hspace{-2em} %to center text under image
	\include{figure5a.tikz}
		\caption{Geometric motivation for the fence.}
		\label{geometric-motivation}
	\end{subfigure}
    \kern10em
	\begin{subfigure}[b]{0.3\textwidth}
    \centering
		\definecolor{p1}{rgb}{0,0,0}
\definecolor{p2}{rgb}{0.1803921568627451,0.49019607843137253,0.19607843137254902}
\definecolor{lcolor}{rgb}{0,.2,.6}
\scalebox{1}{
\begin{tikzpicture}[line cap=round,line join=round,>=triangle 45,x=1cm,y=1.2cm]
\draw [{[-},line width=1.5pt, dotted,color=lcolor] (0,0)-- (0,0.5773502691896257);
\draw [{(-]},line width=1.5pt, color=lcolor] (0,0.5773502691896257)-- (0,2.8867513459481287);
\draw [-),line width=1.5pt, dotted,color=lcolor] (0,2.8867513459481287)-- (0,3.4641016151377544);
\draw (0.15070889992263398,2.0207259421636903) node[anchor=west] {$c_0$};
\draw (0.15070889992263398,3.175426) node[anchor=west] {$c_1$};
\draw (0.15070889992263398,0.8660254) node[anchor=west] {$c_2$};

\begin{scriptsize}
\draw [fill=p1] (0,2.0207259421636903) circle (2.5pt);
\draw [fill=p1] (0,3.175426) circle (2.5pt);
\draw [fill=p1] (0,0.8660254) circle (2.5pt);
\end{scriptsize}
\end{tikzpicture}
}
		%\caption{The fence for the unfolding in (a); the first contact point is labeled $1$.}
                \caption{The fence for the unfolding in (a); contact points  along the unfolding are indexed sequentially.}
		\label{fence}
	\end{subfigure}
	
	\caption{}
    \label{figure4}
\end{figure}

Define $t_i:= i-a$, where $i \in \mathbb Z$; then $c_i:= f(\upsilon\left(  t_{i}\right)) =f ( i,\frac{y}{x}(i-a) )$. Extending the domain of $\upsilon$ to all real numbers allows us to define $c_i$ for all integers $i$, in which case the equality still holds. The following lemma characterizes the purely geometric notion of alignment in terms of analytic conditions on the set of contact points.

\begin{Lemma}\label{align-contact}
The initial and terminal points of an unfolding $\upsilon\left( t\right), \, 0<t \leq T$, are aligned if and only if $c_i =c_{i+T}$ for all $i \in \mathbb Z$.
%$f \left( T, \frac{y}{x} (T-a) \right) = f \left( 0, -\frac{y}{x} a \right)$.
 \end{Lemma}

\begin{proof}
By inspection, the initial point $\left(a, 0\right)$ and the terminal point $\left(T+a, \frac{y}{x}T\right)$ are aligned in the tessellation $\mathcal T$ if and only if  the horizontal change $T$ and the vertical change $ \frac{y}{x}T$ are both integers, and $2\mid \left(T+\frac{y}{x}T\right)$. By definition, $c_i =c_{i+T}$ if and only if $T \in \mathbb Z$ and %$f \left(i, \frac{y}{x} (i-a) \right) = f \left(i+T, \frac{y}{x} (i+T-a) \right)$; the latter can be expressed as 
$i + \frac{y}{x} (i-a) +2 \mathbb Z = i+T +\frac{y}{x}(i+T-a) +2 \mathbb Z$, where the latter condition is equivalent to $2\mid \left(T+\frac{y}{x}T\right)$. The fact that $T \in \mathbb Z$ and $2\mid \left(T+\frac{y}{x}T\right)$ implies $\frac{y}{x}T \in \mathbb Z$ completes the proof.
\end{proof}

Computing the period of a periodic orbit will be significantly simplified by

\begin{Proposition}\label{multiplicity}
    If the initial and terminal points of an unfolding are aligned, all contact points are equally spaced on the fence and have the same multiplicity.
 \end{Proposition}

\begin{proof}
If the initial and terminal points of $\upsilon\left( t\right), \, 0<t \leq T$, are aligned, then $c_i = c_{i+T}$ for all $i \in \mathbb Z$ by Lemma \ref{align-contact}. Hence $\mathcal G := \{ (1+\frac{y}{x})i + 2 \mathbb{Z}:i\in \mathbb Z\}$ is a finite subgroup of $\mathbb{R} / 2 \mathbb{Z}$, and the set of contact points $\{c_i\}=(-\frac{y}{x} a+2 \mathbb{Z})+\mathcal {G}$ is a coset of $\mathcal G$ in $\mathbb R / 2 \mathbb Z$. Therefore contact points are equally spaced in $\mathbb{R}/2\mathbb{Z}$. Furthermore, if $|\mathcal G|=m$, then $m \mid T$ and 
$m_{T}\left(c_i\right)=\frac{T}{m}$ for all $i$.    
\end{proof}

\subsection{Counting the number of edges cut by an unfolding}

Since the initial and terminal points of $\upsilon\left( t\right), \, 0<t \leq T$, are aligned if and only if
both $T$ and  $\frac{y}{x} T$ are integers and $2\mid \left(T+\frac{y}{x}T \right)$, it follows that $x \mid T$ since $\mathrm{gcd}(x,y)=1$. The \emph{first alignment of} $\upsilon$, denoted by $T_1$, is the smallest value of $T$ for which the three conditions above hold. With the parities of $x$ and $y$ in mind, it is easy to check that these three conditions imply

\begin{Proposition}\label{firstalignment}
The first alignment of an unfolding $\upsilon$ is given by 
\begin{align*}	
	T_1=
	\begin{cases}  
              x & \textnormal{if }  x\equiv y\operatorname{mod}2\\
		2x & \textnormal{if } x\not\equiv y\operatorname{mod}2.
	\end{cases}
	\end{align*}
\end{Proposition}

Let $N_T$ be the number of edges of $\mathcal T$ cut by $\upsilon(t) , \, 0< t \leq T$. Proposition \ref{firstalignment} implies that the initial and terminal points of $\upsilon\left( t\right),$ $0< t \leq  2x$, are always aligned. Consequently, we can compute $N_{2x}$ by appealing to the regularity of the contact points given by Proposition \ref{multiplicity}.

\begin{Lemma}\label{decomposition}
Let $b_{2x}$ denote the number of contact points of $\upsilon(t) , \, 0< t \leq 2x$, on the barrier. Then %The number of edges of $\mathcal T$ cut by $\upsilon(t) , \, 0< t \leq 2x$, is 
$$N_{2x}=8y+m_{2x} b_{2x} .$$ 
\end{Lemma}

\begin{proof}
Since $0< t \leq 2x$, the unfolding $\upsilon$ cuts $\frac{y}{x} \cdot 2x =2y$ horizontal inclines. Between consecutive horizontal inclines, $\upsilon$ cuts four non-vertical edges. Recall that the points at which $\upsilon$ cuts vertical edges correspond to the contact points on the barrier. Since the initial and terminal points of $\upsilon$ are aligned, all contact points have the same multiplicity $m_{2x}$ by Proposition \ref{multiplicity}. Thus the total number of contact points on the barrier is $m_{2x} b_{2x}$, which is also the total number of vertical edges cut by $\upsilon$. Consequently, the total number of edges cut by  $\upsilon$ is $4 \cdot 2y+m_{2x} b_{2x}=8y+m_{2x} b_{2x} .$
\end{proof}

An explicit formula for $N_{2x}$ follows from explicit formulas for $m_{ 2x}$ and $b_{ 2x}$. 

\begin{Lemma}\label{spacing}
 The multiplicity $m_{2x}$ and the spacing $s$ between consecutive contact points are given by
    \begin{equation*}
    m_{2x}=\left\{
    \begin{array}
        [c]{cc}%
        2, & \textnormal{if }x\equiv y\operatorname{mod}2\\
        1, & \textnormal{if }x\not\equiv y\operatorname{mod}2%
    \end{array}
    \right.  \text{ and \ }s=\left\{
    \begin{array}
        [c]{cc}%
        2/x, & \textnormal{if } x\equiv y\operatorname{mod}2\\ 
        1/x, & \, \textnormal{if } x\not\equiv y\operatorname{mod}2.%
    \end{array}
    \right.
    \end{equation*}
\end{Lemma}

\begin{proof}
At the first alignment, $\mathcal C_{ T_1}$ consists of $T_1$ distinct contact points each with multiplicity $1$. If $x\equiv y \bmod 2$, then $T_1=x$ by Proposition \ref{multiplicity}, so that each contact point of $\mathcal{C}_{2x}$ has multiplicity $2$. Since there are $x$ distinct contact points on the fence $\mathcal{ F}$ of length $2$,  the spacing $s=\frac{2}{x}$. If $x \not \equiv y \bmod 2$, then $T_1=2x$ so that each contact point in $\mathcal{C}_{2x}$ has multiplicity $1$. Consequently, there are $2x$ distinct contact points and the spacing $s=\frac{2}{2x}=\frac{1}{x}$.
\end{proof}

Having established the spacing $s$, let us derive a formula for the number $b_{2x}$ of contact points on the barrier. Let $\left[x\right]$ denote the integer part of $x$.  

\begin{Lemma}\label{perturb}
The number of contact points on the barrier is
\begin{align*}	
	b_{2x}=
	\begin{cases}  
              \frac{4}{3s} & \textnormal{if }  3 \mid x\smallskip\\
		\left[\frac{4}{3s}\right]  \,\, \textnormal{or} \,\, \left[\frac{4}{3s}\right] +1 & \textnormal{if }  3 \nmid x.
	\end{cases}
	\end{align*}
\end{Lemma}

\begin{proof}
Horizontally translating the initial point $P$ uniformly shifts the contact points and preserves the spacing $s$. Since $s \in\{\frac{1}{x},\frac{2}{x}\}$, the length of barrier $\frac{4}{3}$ is an integer multiple of $s$ if and only if $3 \mid x$. 

Suppose $3 \mid x$. Since contact points are equally spaced and the barrier $\left(  \frac{1}{3},\frac{5}{3}\right]$ is half open, uniformly shifting the contact points preserves the number of contact points on the barrier. Thus $b_{2x}= \frac{4}{3s}$.
Suppose $3\nmid x$. Uniformly shifting the contact points alternately increases or decreases $b_{2x}$ by $1$ as contact points enter or leave the barrier. Therefore  $b_{2x}\in\{\left[\frac{4}{3s}\right],\left[\frac{4}{3s}\right] +1\}$. 
\end{proof}

The facts we need to derive an explicit formula for $N_{ 2x}$ are now in place.

\begin{Proposition}\label{count}
	The number $N : =N_{2x}$ is given by the following table:
	{\small
    \begin{center} 
		\begin{tabular}{| c | c | c || l | }
			\hline
			$N \equiv 1,3 \bmod 4$ & $N \equiv 2 \bmod 4$ & $N \equiv 0 \bmod 4$ & \\ \hline \hline
			&  & $8y+\frac{4x}{3}$ & $x \equiv 0 \bmod 3 \textnormal{ and } x \equiv y \bmod 2$ \\ \hline
			&  & $8y+\frac{4x}{3}$ & $x \equiv 0 \bmod 3 \textnormal{ and } x \not\equiv y \bmod 2$ \\ \hline
			&$8y+\frac{4x+2}{3}$ & $8y+\frac{4x-4}{3}$ & $x \equiv 1 \bmod 3 \textnormal{ and } x \equiv y \bmod 2$\\ \hline
			$8y+\frac{4x-1}{3}$ & $8y+\frac{4x+2}{3}$ & &$x \equiv 1 \bmod 3 \textnormal{ and } x \not\equiv y \bmod 2$ \\ \hline
			&$8y+\frac{4x-2}{3}$ & $8y+\frac{4x+4}{3}$ & $x \equiv 2 \bmod 3 \textnormal{ and } x \equiv y \bmod 2$\\ \hline
			$8y+\frac{4x+1}{3}$ & $8y+\frac{4x-2}{3}$ & &$x \equiv 2 \bmod 3 \textnormal{ and } x \not\equiv y \bmod 2$ \\ \hline
		\end{tabular}
	\end{center}
    }
\end{Proposition}

\begin{remark*}
The columns are arranged to accommodate Proposition~\ref{terminate}.
\end{remark*}

\begin{proof}
The formula for $N_{2x}$ follows immediately from Lemmas \ref{decomposition}, \ref{spacing}, and \ref{perturb}. 
\smallskip

\noindent \textbf{Case 1}. If $x \equiv 0 \bmod 3$ and $x \equiv y \bmod 2$, then $m_{2x}=2$, $s=\frac{2}{x}$, and $b_{2x}=\frac{4}{3s}$. Thus $N_{2x}=8y+m_{2x} b_{2x}=8y+ 2 \cdot \frac{4}{3} \cdot \frac{x}{2}=8y+ \frac{4x}{3}$.
\smallskip

\noindent \textbf{Case 4}.  If $x \equiv 1 \bmod 3$ and $x \not\equiv y \bmod 2$, then $m_{2x}=1$, $s=\frac{1}{x}$, and either $b_{2x}=\left[\frac{4}{3s}\right]$ or $\left[\frac{4}{3s}\right] +1$. If $b_{2x}=\left[\frac{4}{3s}\right]$, then $N_{2x}=8y+m_{2x} b_{2x}=8y+ 1 \cdot \left[\frac{4}{3} \cdot \frac{x}{1}\right]=8y+\frac{4x-1}{3}$. If $b_{2x}=\left[\frac{4}{3s}\right]+1$, then $N_{2x}=8y+ 1 \cdot \left( \left[\frac{4}{3}\cdot \frac{x}{1}\right]+1 \right)= 8y+\frac{4x+2}{3}$.
\smallskip

\noindent Proofs of the other cases are similar and left to the reader.	
\end{proof}

\subsection{Computing the period of a periodic orbit}

Since the period of a periodic orbit is the period of its primitive sub-orbits, let us determine the values of $T$ for which $\upsilon (t), 0<t \leq T$, is an unfolding of a primitive periodic orbit.

\begin{Proposition}\label{terminate}
	Let $x,y\in \mathbb{N}$ such that $x < y$ and $\mathrm{gcd}(x,y)=1$. A primitive periodic orbit $\gamma$ with initial angle $\Theta = \arctan(\frac{y}{x}\sqrt{3})$ has an unfolding $\upsilon(t)$, $0 < t \leq T$, for some $T\in\{x,2x,4x\}$, and period $p(x,y)$ determined as follows:
    \begin{enumerate}
      \item If $x\equiv y \bmod 2$ and $N_{2x} \equiv 0 \bmod 4$, then $T=x$ and $p(x,y)=\frac{1}{2} N_{2x}$.
        \item If $x\equiv y \bmod 2$ and $N_{2x}\equiv 2 \bmod 4$, then $T=2x$ and $p(x,y)=N_{2x}$.
        \item If $x \not\equiv y \bmod 2$ and $N_{2x}$ is even, then $T=2x$ and $p(x,y) = N_{2x}$.
        \item If $x\not\equiv y \bmod 2$ and $N_{2x}$ is odd, then $T=4x$ and $p(x,y)=2 N_{2x}$.
     
    \end{enumerate}
\end{Proposition}

\begin{proof}
	A primitive periodic orbit has an unfolding $\upsilon(t)$, $0 < t \leq T$, if and only if  $T$ is the smallest positive integer such that the initial and terminal points are aligned and $N_T$ is even. Recall that the alignment condition implies $x \mid T$, and the first alignment $T_1=x$ when $x\equiv y \bmod 2$ and $T_1=2x$ otherwise.		
\smallskip

	\noindent\textbf{Case 1}. If $x\equiv y \bmod 2$, then $T_1=x$. If $N_{2x} \equiv 0 \bmod 4$, then $N_x=\frac{1}{2} N_{2x}$ is even so that $T=x$ and $p(x,y)=N_x=\frac{1}{2}N_{2x}$. If $N_{2x}\equiv 2 \bmod 4$, then $N_x=\frac{1}{2} N_{2x}$ is odd so that $T=2x$ and $p(x,y)=N_{2x}$.

\smallskip
	\noindent\textbf{Case 2}. If $x\not\equiv y \bmod 2$, then $T_1=2x$. If $N_{2x}$ is even, then $T=2x$ and $p(x,y)=N_{2x}$. If $N_{2x}$ is odd, then $N_{4x}=2 N_{2x}$ is even so that  $T=4x$ and $p(x,y)=N_{4x}=2N_{2x}$.
\end{proof}

The period of every periodic orbit is now determined.

\begin{Theorem}\label{formula}
	Let $x,y \in \mathbb{N}$ such that $x < y$ and $\mathrm{gcd}(x,y)=1$. Then the period $p(x,y)$ of a periodic orbit with initial angle $\Theta=\arctan(\frac{y}{x}\sqrt{3})$ is 
    $$p(x, y)=
	\begin{cases}
	4y + \frac{2x}{3} &  {\rm if} \,\, x \equiv 0  \textnormal{ mod 3}  \textnormal{ and }  x \equiv y \textnormal{ mod 2}\\
	8y + \frac{4x}{3} &  {\rm if} \,\, x \equiv 0  \textnormal{ mod 3}  \textnormal{ and }  x \not\equiv y \textnormal{ mod 2}\\
	4y + \frac{2x-2}{3} \textnormal{ or } \textnormal 8y + \frac{4x+2}{3} & {\rm if} \,\, x \equiv 1  \textnormal{ mod 3}  \textnormal{ and }  x \equiv y \textnormal{ mod 2}\\
	16y + \frac{8x-2}{3} \textnormal{ or } 8y + \frac{4x+2}{3} & {\rm if} \,\, x \equiv 1  \textnormal{ mod 3}  \textnormal{ and }  x \not\equiv y \textnormal{ mod 2}\\
	4y + \frac{2x+2}{3} \textnormal{ or } 8y + \frac{4x-2}{3} & {\rm if} \,\, x \equiv 2  \textnormal{ mod 3}  \textnormal{ and }  x \equiv y \textnormal{ mod 2}\\
	16y + \frac{8x+2}{3} \textnormal{ or } 8y + \frac{4x-2}{3} & {\rm if} \,\, x \equiv 2  \textnormal{ mod 3}  \textnormal{ and }  x \not\equiv y \textnormal{ mod 2}.\\
	\end{cases}$$
\end{Theorem}

\begin{proof}
Propositions \ref{count} and \ref{terminate} immediately lead to the formula for $p(x,y)$. 

\smallskip

\noindent \textbf{Case 1}. If $x \equiv 0 \bmod 3$ and $x \equiv y \bmod 2$, then $N_{2x}=8y+\frac{4x}{3} \equiv 0 \bmod 4$ so that $p(x,y)=\frac{1}{2} N_{2x}= 4y+\frac{2x}{3}$.
\smallskip

\noindent \textbf{Case 4}.  If $x \equiv 1 \bmod 3$ and $x \not\equiv y \bmod 2$, then $N_{2x}=8y+\frac{4x-1}{3}$ or  $N_{2x}=8y+\frac{4x+2}{3}$. In the first case, $N_{2x}$ is odd so that $p(x,y)=2 N_{2x}= 16y+\frac{8x-2}{3}$; in the second case, $N_{2x}$ is even so that $p(x,y)=N_{2x}= 8y+\frac{4x+2}{3}$.
\smallskip

\noindent Proofs of the other cases are similar and left to the reader.	
\end{proof}

\noindent When $\Theta=60^\circ$, periods $4$ and $10$ are consistent with the formula in Theorem~\ref{formula} with $(x,y)=(1,1)$; when $\Theta=90^\circ$, period $8$ is consistent with $(x,y)=(0,1)$.
	
\begin{Corollary}
	A biperiodic periodic orbit with initial angle $\Theta = \arctan \big(\frac{y}{x} \sqrt{3} \big)$ has one of two possible periods $p_1 < p_2$, where $p_2  = 2p_1 + 2$ or $p_2 = 2 p_1 - 2$.
\end{Corollary}

While our understanding of unfoldings follows by considering the fence, the simple statement in Thoerem 12 
makes no mention of the initial point.  However, the framework developed here allows us to determine a more precise period
formula in terms of the initial and terminal points $P = (a, 0)$ and $Q = (a + x, y)$, and doing so  required us to compute the number of contact points on the barrier as a function of $a$:
 $$b_a(x,y)=\left[ \frac{5/3 + a y/x}{s_a(x,y)} \right]-\left[ \frac{1/3 + a y/x}{s_a(x,y)} \right],$$
  where $s_a(x,y)$ is the spacing function.

\section{Periodic orbits on other obtuse polygons}\label{othershapes}

The methods developed in Section \ref{120isosceles} can be applied to a $60^\circ$-rhombus and a $60^\circ$-$90^\circ$-$120^\circ$-kite.  Let $\mathcal{T}_1$ and $\mathcal{T}_2$ be the respective edge tessellations generated by a $60^\circ$-rhombus and a $60^\circ$-$90^\circ$-$120^\circ$-kite. Note that an analogue of Theorem \ref{align+orient} holds in both of these cases. In either case, periodic orbits can be represented by unfoldings with an initial angle $\Theta \in [60^\circ, 90^\circ]$, and $\Theta \in (60^\circ, 90^\circ)$ can be expressed in the form $\Theta=\arctan(\frac{y}{x}\sqrt{3})$ with $x < y$ and $\gcd(x, y)=1$. When $\Theta=60^\circ$ or $\Theta=90^\circ$, the period can be determined by inspection and fits the general formula to be derived.

\subsection{The $60^\circ$-rhombus}

Note that $\mathcal {T}_1$ can be obtained from $\mathcal T$ by removing its $0^\circ$, $60^\circ$, and $120^\circ$ inclines (see Figure~\ref{60-rhombusT}). Impose the same coordinate system on $\mathcal{T}_1$ we imposed on $\mathcal{T}$. The barrier and gate for $\mathcal{T}_1$ are identical to the barrier and gate for $\mathcal{T}$, and all definitions and techniques in the previous sections apply. Although  $\mathcal{T}_1$ has no horizontal inclines, we can position the initial point of an unfolding $\upsilon$ on a horizontal incline of $\mathcal T$ (the dotted line on Figure \ref{60-rhombusT}). Although a $60^\circ$-rhombus exhibits both line and rotational symmetry, a quick check shows that the initial and terminal rhombuses determined by an unfolding with aligned initial and terminal points, may differ by a reflection but not by a rotation.

%\begin{figure}[ht!]
  %  \centering
  %  \include{rhombus.tikz}
   % \caption{The tessellation $\mathcal T_1$ generated by a $60^\circ$-rhombus. The dotted lines are horizontal inclines of $\mathcal T$ but not $\mathcal T_1$.}
   % \label{60-rhombusT}
%\end{figure}

The formula $N_{2x}=4y+m_{2x} b_{2x}$ (as in Lemma~\ref{decomposition}) leads to the following table (as in Proposition~\ref{count}):
{\small \begin{center}
		\begin{tabular}{| c | c | c || l | } \hline
			$N \equiv 1,3 \bmod 4$ & $N \equiv 2 \bmod 4$ & $N \equiv 0 \bmod 4$ & \\ \hline \hline
			&  & $4y+\frac{4x}{3}$ & $x \equiv 0 \bmod 3 \textnormal{ and } x \equiv y \bmod 2$ \\ \hline
			&  & $4y+\frac{4x}{3}$ & $x \equiv 0 \bmod 3 \textnormal{ and } x \not\equiv y \bmod 2$ \\ \hline
			&$4y+\frac{4x+2}{3}$& $4y+\frac{4x-4}{3}$ & $x \equiv 1 \bmod 3 \textnormal{ and } x \equiv y \bmod 2$\\ \hline
			$4y+\frac{4x-1}{3}$ & $4y+\frac{4x+2}{3}$ & &$x \equiv 1 \bmod 3 \textnormal{ and } x \not\equiv y \bmod 2$ \\ \hline
			&$4y+\frac{4x-2}{3}$& $4y+\frac{4x+4}{3}$ & $x \equiv 2 \bmod 3 \textnormal{ and } x \equiv y \bmod 2$\\ \hline
			$4y+\frac{4x+1}{3}$ & $4y+\frac{4x-2}{3}$ & &$x \equiv 2 \bmod 3 \textnormal{ and } x \not\equiv y \bmod 2$ \\ \hline
		\end{tabular}
	\end{center}
}

Proposition~\ref{terminate} applies, and when combined with the table above, produces the following formula for the period (as in Theorem~\ref{formula}):	
$$p(x, y)=
	\begin{cases}
	2y + \frac{2x}{3} &  {\rm if} \,\, x \equiv 0  \textnormal{ mod 3}  \textnormal{ and }  x \equiv y \textnormal{ mod 2}\\
	4y + \frac{4x}{3} &  {\rm if} \,\, x \equiv 0  \textnormal{ mod 3}  \textnormal{ and }  x \not\equiv y \textnormal{ mod 2}\\
	2y + \frac{2x-2}{3} \textnormal{ or } \textnormal 4y + \frac{4x+2}{3} & {\rm if} \,\, x \equiv 1  \textnormal{ mod 3}  \textnormal{ and }  x \equiv y \textnormal{ mod 2}\\
	4y + \frac{4x+2}{3} \textnormal{ or } 8y + \frac{8x-2}{3} & {\rm if} \,\, x \equiv 1  \textnormal{ mod 3}  \textnormal{ and }  x \not\equiv y \textnormal{ mod 2}\\
	2y + \frac{2x+2}{3} \textnormal{ or } 4y + \frac{4x-2}{3} & {\rm if} \,\, x \equiv 2  \textnormal{ mod 3}  \textnormal{ and }  x \equiv y \textnormal{ mod 2}\\
	4y + \frac{4x-2}{3} \textnormal{ or }8y + \frac{8x+2}{3} & {\rm if} \,\, x \equiv 2  \textnormal{ mod 3}  \textnormal{ and }  x \not\equiv y \textnormal{ mod 2}.\\
	\end{cases}$$

\begin{figure}[]
	\centering
	
	\begin{subfigure}[b]{0.4\textwidth}
    \centering
    %\hspace{-4em} %to center text under image
	\definecolor{niceblue}{rgb}{0,.2,.6}
\scalebox{.75}{%
\begin{tikzpicture}[line cap=round,line join=round,>=triangle 45,x=.9cm,y=.9cm]
\draw [color=niceblue] (0,0)-- (1,0.5773502691896257);
\draw [color=niceblue] (1,1.7320508075688772)-- (0,1.1547005383792515);
\draw [color=niceblue] (0,1.1547005383792515)-- (0,0);
\draw [color=niceblue] (1,0.5773502691896257)-- (1,1.7320508075688772);
\draw [color=niceblue] (1,1.7320508075688772)-- (2,1.1547005383792512);
\draw [color=niceblue] (2,0)-- (3,0.5773502691896258);
\draw [color=niceblue] (3,1.7320508075688774)-- (2,1.1547005383792512);
\draw [color=niceblue] (2,1.1547005383792512)-- (2,0);
\draw [color=niceblue] (6,0)-- (5,0.5773502691896273);
\draw [color=niceblue] (5,1.7320508075688787)-- (6,1.1547005383792537);
\draw [color=niceblue] (6,1.1547005383792537)-- (6,0);
\draw [color=niceblue] (4,0)-- (5,0.5773502691896273);
\draw [color=niceblue] (5,0.5773502691896273)-- (5,1.7320508075688787);
\draw [color=niceblue] (5,1.7320508075688787)-- (4,1.154700538379252);
\draw [color=niceblue] (4,0)-- (3,0.5773502691896258);
\draw [color=niceblue] (3,0.5773502691896258)-- (3,1.7320508075688774);
\draw [color=niceblue] (3,1.7320508075688774)-- (4,1.154700538379252);
\draw [color=niceblue] (4,1.154700538379252)-- (4,0);
\draw [color=niceblue] (1,1.732050807568876)-- (1,0.5773502691896246);
\draw [color=niceblue] (1,0.5773502691896246)-- (2,0);
\draw [dotted,color=niceblue] (0,0)-- (6,0);
\draw [dotted,color=niceblue] (0,1.7320508075688772)-- (6,1.7320508075688772);
\draw [color=niceblue] (0,3.464101615137755)-- (1,2.886751345948129);
\draw [color=niceblue] (1,1.7320508075688776)-- (0,2.3094010767585034);
\draw [color=niceblue] (0,2.3094010767585034)-- (0,3.464101615137755);
\draw [color=niceblue] (1,2.886751345948129)-- (1,1.7320508075688776);
\draw [color=niceblue] (1,1.7320508075688776)-- (2,2.309401076758504);
\draw [color=niceblue] (2,3.4641016151377553)-- (3,2.886751345948129);
\draw [color=niceblue] (3,1.7320508075688774)-- (2,2.309401076758504);
\draw [color=niceblue] (2,2.309401076758504)-- (2,3.4641016151377553);
\draw [color=niceblue] (6,3.4641016151377526)-- (5,2.8867513459481273);
\draw [color=niceblue] (5,1.732050807568876)-- (6,2.309401076758501);
\draw [color=niceblue] (6,2.309401076758501)-- (6,3.4641016151377526);
\draw [color=niceblue] (4,3.4641016151377544)-- (5,2.8867513459481273);
\draw [color=niceblue] (5,2.8867513459481273)-- (5,1.732050807568876);
\draw [color=niceblue] (5,1.732050807568876)-- (4,2.309401076758503);
\draw [color=niceblue] (4,3.4641016151377544)-- (3,2.886751345948129);
\draw [color=niceblue] (3,2.886751345948129)-- (3,1.7320508075688774);
\draw [color=niceblue] (3,1.7320508075688774)-- (4,2.309401076758503);
\draw [color=niceblue] (4,2.309401076758503)-- (4,3.4641016151377544);
\draw [color=niceblue] (1,1.7320508075688787)-- (1,2.88675134594813);
\draw [color=niceblue] (1,2.88675134594813)-- (2,3.4641016151377553);
\draw [dotted,color=niceblue] (0,3.464101615137755)-- (6,3.464101615137755);
\draw [color=niceblue] (0,3.464101615137755)-- (1,4.041451884327381);
\draw [color=niceblue] (1,5.196152422706632)-- (0,4.618802153517006);
\draw [color=niceblue] (0,4.618802153517006)-- (0,3.464101615137755);
\draw [color=niceblue] (1,4.041451884327381)-- (1,5.196152422706632);
\draw [color=niceblue] (1,5.196152422706632)-- (2,4.618802153517006);
\draw [color=niceblue] (2,3.4641016151377544)-- (3,4.041451884327381);
\draw [color=niceblue] (3,5.196152422706632)-- (2,4.618802153517006);
\draw [color=niceblue] (6,3.464101615137757)-- (5,4.041451884327382);
\draw [color=niceblue] (5,5.196152422706634)-- (6,4.6188021535170085);
\draw [color=niceblue] (6,4.6188021535170085)-- (6,3.464101615137757);
\draw [color=niceblue] (4,3.4641016151377553)-- (5,4.041451884327382);
\draw [color=niceblue] (5,4.041451884327382)-- (5,5.196152422706634);
\draw [color=niceblue] (5,5.196152422706634)-- (4,4.618802153517007);
\draw [color=niceblue] (4,3.4641016151377553)-- (3,4.041451884327381);
\draw [color=niceblue] (3,4.041451884327381)-- (3,5.196152422706632);
\draw [color=niceblue] (3,5.196152422706632)-- (4,4.618802153517007);
\draw [color=niceblue] (4,4.618802153517007)-- (4,3.4641016151377553);
\draw [color=niceblue] (1,5.196152422706631)-- (1,4.04145188432738);
\draw [color=niceblue] (1,4.04145188432738)-- (2,3.4641016151377544);
\draw [color=niceblue] (2,3.4641016151377544)-- (2,4.618802153517006);
\draw [dotted,color=niceblue] (0,5.196152422706632)-- (6,5.196152422706632);
\end{tikzpicture}
}
    \caption{The tessellation $\mathcal T_1$ generated by a $60^\circ$-rhombus.}
    \label{60-rhombusT}
	\end{subfigure}
    \kern2em
	\begin{subfigure}[b]{0.4\textwidth}
    \centering
    \definecolor{niceblue}{rgb}{0,0.2,0.6}
\begin{tikzpicture}[scale=.5][line cap=round,line join=round,>=triangle 45,x=.7cm,y=.7cm]
\clip (0,0) rectangle  (8,6.928203230275522);
\draw [color=niceblue] (0,0)-- (2,0);
\draw [color=niceblue] (2,0)-- (2,1.1547005383792515);
\draw [color=niceblue] (2,1.1547005383792515)-- (1,1.7320508075688772);
\draw [color=niceblue] (1,1.7320508075688772)-- (0,0);
\draw [color=niceblue] (4,0)-- (2,0);
\draw [color=niceblue] (2,1.1547005383792515)-- (3,1.732050807568877);
\draw [color=niceblue] (3,1.732050807568877)-- (4,0);
\draw [color=niceblue] (3,1.7320508075688772)-- (2,3.464101615137755);
\draw [color=niceblue] (1,1.7320508075688772)-- (2,3.4641016151377544);
\draw [color=niceblue] (6,3.464101615137751)-- (5,1.7320508075688752);
\draw [color=niceblue] (5,1.7320508075688752)-- (4,2.3094010767585025);
\draw [color=niceblue] (4,2.3094010767585025)-- (4,3.4641016151377535);
\draw [color=niceblue] (4,3.4641016151377535)-- (6,3.464101615137751);
\draw [color=niceblue] (4,2.3094010767585025)-- (3,1.7320508075688776);
\draw [color=niceblue] (4,3.4641016151377535)-- (2,3.4641016151377566);
\draw [color=niceblue] (4,0)-- (5,1.7320508075688774);
\draw [color=niceblue] (6,1.154700538379246)-- (7,1.732050807568869);
\draw [color=niceblue] (7,1.732050807568869)-- (6,3.464101615137749);
\draw [color=niceblue] (6,0)-- (4,0);
\draw [color=niceblue] (6,1.1547005383792532)-- (6,0);
\draw [color=niceblue] (5,1.7320508075688772)-- (6,1.1547005383792532);
\draw [color=niceblue] (6,0)-- (8,0);
\draw [color=niceblue] (8,0)-- (7,1.7320508075688807);
\draw [color=niceblue] (0,2.309401076758503)-- (1,1.7320508075688772);
\draw [color=niceblue] (0,2.309401076758503)-- (0,3.4641016151377544);
\draw [color=niceblue] (0,3.4641016151377544)-- (2,3.4641016151377557);
\draw [color=niceblue] (0,4.618802153517006)-- (1,5.196152422706633);
\draw [color=niceblue] (1,5.196152422706633)-- (2,3.464101615137757);
\draw [color=niceblue] (0,4.618802153517006)-- (0,3.4641016151377544);
\draw [color=niceblue] (2,5.77350269189626)-- (1,5.196152422706634);
\draw [color=niceblue] (2,5.77350269189626)-- (3,5.1961524227066365);
\draw [color=niceblue] (3,5.196152422706636)-- (2,3.4641016151377566);
\draw [color=niceblue] (4,4.618802153517009)-- (4,3.464101615137758);
\draw [color=niceblue] (4,4.618802153517009)-- (3,5.196152422706634);
\draw [color=niceblue] (5,5.196152422706639)-- (6,3.4641016151377633);
\draw [color=niceblue] (5,5.196152422706641)-- (4,6.928203230275514);
\draw [color=niceblue] (4,4.618802153517009)-- (5,5.196152422706637);
\draw [color=niceblue] (7,5.196152422706644)-- (6,3.4641016151377624);
\draw [color=niceblue] (8,2.30940107675849)-- (8,3.4641016151377415);
\draw [color=niceblue] (8,4.6188021535169925)-- (8,3.4641016151377415);
\draw [color=niceblue] (8,2.30940107675849)-- (7,1.7320508075688674);
\draw [color=niceblue] (8,4.618802153517019)-- (7,5.196152422706643);
\draw [color=niceblue] (6,5.773502691896266)-- (7,5.196152422706644);
\draw [color=niceblue] (6,5.773502691896266)-- (5,5.19615242270664);
\draw [color=niceblue] (8,3.4641016151377415)-- (6,3.464101615137749);
\draw [color=niceblue] (4,6.928203230275509)-- (2,6.928203230275514);
\draw [color=niceblue] (0,6.928203230275509)-- (2,6.928203230275511);
\draw [color=niceblue] (1,5.196152422706633)-- (0,6.928203230275509);
\draw [color=niceblue] (3,5.196152422706638)-- (4,6.928203230275516);
\draw [color=niceblue] (2,6.928203230275511)-- (2,5.77350269189626);
\draw [color=niceblue] (4,6.928203230275509)-- (6,6.92820323027552);
\draw [color=niceblue] (6,6.928203230275518)-- (6,5.773502691896266);
\draw [color=niceblue] (7,5.196152422706643)-- (8,6.928203230275522);
\draw [color=niceblue] (8,6.928203230275522)-- (6,6.928203230275518);

\begin{scriptsize}

\draw (0,0) node[anchor= south west]{$A$};
\draw (2.1,0) node[anchor= south west]{$O$};
\draw (2,1.8*0.7) node[anchor=south]{$B$};
\draw (1,2.8*0.5773502691896257) node[anchor=east]{$D$};

\end{scriptsize}

\end{tikzpicture}
    
    \caption{The tessellation $\mathcal T_2$ generated by a $60^\circ$-$90^\circ$-$120^\circ$-kite.}
    \label{60-90-120-kiteT}
	\end{subfigure}
	
	\caption{}
    \label{figure5}
\end{figure}

\subsection{The $60^\circ$-$90^\circ$-$120^\circ$-kite}

The edge tessellation $\mathcal{T}_2$ is also related to $\mathcal T$. Consider the $60^\circ$-$90^\circ$-$120^\circ$-kite positioned as $\Box AOBD$, where $A$ and $B$ coincide with the two vertices of the $120^\circ$-isosceles $\triangle ABC$ and $O$ is the midpoint of $\overline{AC}$  (see Figure~\ref{60-90-120-kiteT}). Impose the same coordinate system on $\mathcal{T}_2$ we imposed on $\mathcal{T}$ and position the initial point of an unfolding $\upsilon$ on a horizontal incline. While the fence is the same as before, the barrier and gate are interchanged, i.e., $(\frac{1}{3},\frac{5}{3}]$ is the gate.

%\begin{figure}[ht!]
    %\centering
    %\include{kite.tikz}
    %\caption{The tessellation $\mathcal T_2$ generated by a $60^\circ$-$90^\circ$-$120^\circ$-kite.}
   % \label{60-90-120-kiteT}
%\end{figure}

The formula $N_{2x}=6y+m_{2x} b_{2x}$ (as in Lemma~\ref{decomposition}) leads to the following table (as in Proposition~\ref{count}):

{\small \begin{center} 
		\begin{tabular}{| c | c | c || l | } \hline
			$N \equiv 1,3 \bmod 4$ & $N \equiv 2 \bmod 4$ & $N \equiv 0 \bmod 4$ & \\ \hline \hline
			&  & $6y+\frac{2x}{3}$ & $x \equiv 0 \bmod 3 \textnormal{ and } x \equiv y \bmod 2$ \\ \hline
			& $6y+\frac{2x}{3}$ &  & $x \equiv 0 \bmod 3 \textnormal{ and } x \not\equiv y \bmod 2$ \\ \hline
			& $6y+\frac{2x-2}{3}$& $6y+\frac{2x+4}{3}$ & $x \equiv 1 \bmod 3 \textnormal{ and } x \equiv y \bmod 2$\\ \hline
			$6y+\frac{2x+1}{3}$ & &$6y+\frac{2x-2}{3}$  &$x \equiv 1 \bmod 3 \textnormal{ and } x \not\equiv y \bmod 2$ \\ \hline
			&$6y+\frac{2x+2}{3}$& $6y+\frac{2x-4}{3}$ & $x \equiv 2 \bmod 3 \textnormal{ and } x \equiv y \bmod 2$\\ \hline
			$6y+\frac{2x-1}{3}$ & & $6y+\frac{2x+2}{3}$  &$x \equiv 2 \bmod 3 \textnormal{ and } x \not\equiv y \bmod 2$ \\ \hline
		\end{tabular}
	\end{center}
    }

Proposition~\ref{terminate} applies, and when combined with the table above, produces the following formula for the period (as in Theorem~\ref{formula}):	 	
$$p(x, y)=
	\begin{cases}
	3y + \frac{x}{3} &  {\rm if} \,\, x \equiv 0  \textnormal{ mod 3}  \textnormal{ and }  x \equiv y \textnormal{ mod 2}\\
	6y + \frac{2x}{3} &  {\rm if} \,\, x \equiv 0  \textnormal{ mod 3}  \textnormal{ and }  x \not\equiv y \textnormal{ mod 2}\\
	3y + \frac{x+2}{3} \textnormal{ or } \textnormal 6y + \frac{2x-2}{3} & {\rm if} \,\, x \equiv 1  \textnormal{ mod 3}  \textnormal{ and }  x \equiv y \textnormal{ mod 2}\\
	6y + \frac{2x-2}{3} \textnormal{ or } 12y + \frac{4x+2}{3} & {\rm if} \,\, x \equiv 1  \textnormal{ mod 3}  \textnormal{ and }  x \not\equiv y \textnormal{ mod 2}\\
	3y + \frac{x-2}{3} \textnormal{ or } 6y + \frac{2x+2}{3} & {\rm if} \,\, x \equiv 2  \textnormal{ mod 3}  \textnormal{ and }  x \equiv y \textnormal{ mod 2}\\
	6y + \frac{2x+2}{3} \textnormal{ or }12y + \frac{4x-2}{3} & {\rm if} \,\, x \equiv 2  \textnormal{ mod 3}  \textnormal{ and }  x \not\equiv y \textnormal{ mod 2}.\\
	\end{cases}$$

\subsection{The Regular Hexagon}

While the techniques in Section \ref{120isosceles} can be applied to determine the number of edges cut by an unfolding of a periodic orbit on a regular hexagon, the presence of rotational symmetries renders this information insufficient to determine exactly where an unfolding terminates. Position the hexagon so that a subdivision of its edge tessellation is the edge tessellation generated by a $120^\circ$-isosceles triangle, and restrict considerations to initial angles $\Theta\in(30^\circ,60^\circ)$.

Although we are unable to rigorously solve the problem, we pose a conjecture based on extensive numerical evidence. We computed the period for 3814 pairs $(x,y)$ and partitioned these pairs into planar groups using a mixture algorithm \cite{G}. We conjecture that the period $p(x,y)$ is
\begin{equation*}
	\begin{cases}
	3y + x &  {\rm if} \,\, 
	x \equiv 0  \textnormal{ mod 3}, x \equiv y \textnormal{ mod 2}, \textnormal{and } (x,y) \in A_{0,0} \\
	y + \frac{x}{3} &  {\rm if} \,\, 
	x \equiv 0  \textnormal{ mod 3}, x \equiv y \textnormal{ mod 2}, \textnormal{and } (x,y) \in A^c_{0,0}  \\
	6y + 2x &  {\rm if} \,\, 
	x \equiv 0  \textnormal{ mod 3}, x \not\equiv y \textnormal{ mod 2}, \textnormal{and } (x,y) \in A_{0,1} \\
	2y + \frac{2x}{3} &  {\rm if} \,\, 
	x \equiv 0  \textnormal{ mod 3}, x \not\equiv y \textnormal{ mod 2}, \textnormal{and } (x,y) \in A^c_{0,1} \\
	2y + \frac{2x-2}{3} \textnormal{ or } \textnormal 3y + x+2 & {\rm if} \,\, 
	x \equiv 1  \textnormal{ mod 3}, x \equiv y \textnormal{ mod 2}, \textnormal{and } (x,y) \in A_{1,0} \\
	y + \frac{x+2}{3} \textnormal{ or } \textnormal 2y + \frac{2x-2}{3} & {\rm if} \,\, 
	x \equiv 1  \textnormal{ mod 3}, x \equiv y \textnormal{ mod 2}, \textnormal{and } (x,y) \in A^c_{1,0} \\
	4y + \frac{4x+2}{3} \textnormal{ or } 6y + 2x-2 & {\rm if} \,\, 
	x \equiv 1  \textnormal{ mod 3}, x \not\equiv y \textnormal{ mod 2}, \textnormal{and } (x,y) \in A_{1,1} \\
	2y + \frac{2x-2}{3} \textnormal{ or } 4y + \frac{4x+2}{3} & {\rm if} \,\,
	x \equiv 1  \textnormal{ mod 3}, x \not\equiv y \textnormal{ mod 2}, \textnormal{and } (x,y) \in A^c_{1,1} \\
	2y + \frac{2x+2}{3} \textnormal{ or } 3y + x-2 & {\rm if} \,\, 
	x \equiv 2  \textnormal{ mod 3}, x \equiv y \textnormal{ mod 2}, \textnormal{and } (x,y) \in A_{2,0} \\
	y + \frac{x-2}{3} \textnormal{ or } 2y + \frac{2x+2}{3} & {\rm if} \,\, 
	x \equiv 2  \textnormal{ mod 3}, x \equiv y \textnormal{ mod 2}, \textnormal{and } (x,y) \in A^c_{2,0} \\
	4y + \frac{4x-2}{3} \textnormal{ or } 6y + 2x+2 & {\rm if} \,\, 
	x \equiv 2  \textnormal{ mod 3}, x \not\equiv y \textnormal{ mod 2}, \textnormal{and } (x,y) \in A_{2,1} \\
	2y + \frac{2x+2}{3} \textnormal{ or } 4y + \frac{4x-2}{3} & {\rm if} \,\, 
	x \equiv 2  \textnormal{ mod 3}, x \not\equiv y \textnormal{ mod 2}, \textnormal{and } (x,y) \in A^c_{2,1},
	\end{cases}
\end{equation*}
where the sets $A_{i,j}$ have yet to be determined. A numerical grid search appeared to indicate that the sets $A_{i,j}$ cannot be described by a linear modulus condition $c_1 x + c_2 y \bmod c_3$ for any $c_1,c_2=-36,-35,...,35,36$ and $c_3=2,...,36$. However, when $3 \mid x$, we conjecture that whenever $(x,y) \in A_{i,j}$ for any $i,j$, then $(27y-7x,11y-3x)\in A_{i',j'}$ for some possibly different $i',j'$. Under this assumption, we can state our conjecture for the period formula more precisely for a given $x$. %after a small number of periods are found manually.

\section*{Acknowledgments}

We wish to thank David Brown for his participation in an undergraduate research seminar during the 2011-2012 academic year in which this problem was first considered, and Joshua Pavoncello for writing a computer program \cite{P} that visualizes the orbits on an edge tessellating polygon and collects related experimental data. And we wish to thank the referees whose numerous helpful suggestions improved the exposition.

\section*{About the authors}

\subsection*{Benjamin R. Baer}
Department of Statistics and Data Science, Cornell University, Ithaca, NY 14853. brb225@cornell.edu

Mr. Baer is a sixth year statistics PhD student at Cornell University studying statistical methodology.

\subsection*{Faheem Gilani}~
Department of Mathematics, The Pennsylvania State University, State College, PA 16801. fhg3@psu.edu

Mr. Gilani is a fifth year mathematics PhD student at Penn State studying machine learning. 

\subsection*{Zhigang Han}
Department of Mathematics, Millersville University of Pennsylvania, Millersville, PA 17551. zhigang.han@millersville.edu

Dr. Han is an Associate Professor of Mathematics at Millersville University. His primary research areas are symplectic geometry and topology.

\subsection*{Ronald Umble}
Department of Mathematics, Millersville University of Pennsylvania, Millersville, PA 17551. ron.umble@millersville.edu

Dr. Umble retired from teaching in August of 2020. He is a Professor of Mathematics Emeritus at Millersville University and a member of Pi Mu Epsilon–Pennsylvania Zeta Chapter.  His primary research area is algebraic topology.  He has directed or codirected 19 undergraduate research projects, nine of which led to publications. 

\end{document}